 \documentclass[final,3p,times]{elsarticle}

\usepackage{amssymb}
 \usepackage{amsthm, amsmath}
 \usepackage[T1]{fontenc}
 
\usepackage[ansinew]{inputenc} 
\usepackage{mathrsfs}
\usepackage{euscript}
\usepackage{natbib}

\usepackage[draft]{showlabels}
 
\newtheorem{theorem}{Theorem}
\newtheorem{definition}{Definition}
\newtheorem{proposition}{Proposition}
\newdefinition{rmk}{ Remark}
\newproof{pf}{Proof of Theorem 1}

\begin{document}

\begin{frontmatter}



\title{ Blow-up of solutions to semilinear wave equations with a   time-dependent  strong damping}


 \author[Fino]{Ahmad Z. Fino}
 \address[Fino]{Department of Mathematics, Sultan Qaboos University,  FracDiff Research Group (DR/RG/03),  P.O. Box 46, Al-Khoud 123, Muscat, Oman}


 \ead{ahmad.fino01@gmail.com; a.fino@squ.edu.om}
 
  \author[Hamza]{Mohamed Ali Hamza}
\address[Hamza]{Basic Sciences Department, Deanship of Preparatory Year and Supporting Studies, P. O. Box 1982, Imam Abdulrahman Bin Faisal University, Dammam, KSA}
 
  \ead{mahamza@iau.edu.sa}
  
 \begin{abstract}
The paper investigates a class of a semilinear wave equation with time-dependent damping  term ($-\frac{1}{{(1+t)}^{\beta}}\Delta u_t$) and a nonlinearity $|u|^p$. We will show 
the influence of the  the parameter  $\beta$  in the blow-up results under some hypothesis on the initial data and the exponent $p$ by using the test function method. We also study the local existence in time of mild solution in the energy space $H^1(\mathbb{R}^n)\times L^2(\mathbb{R}^n)$.
\end{abstract}


\end{frontmatter}


\section{Introduction}
\setcounter{equation}{0} 
The aim of the paper is to establish a blow-up result for local in time solutions to the Cauchy
 problem for the following semilinear strong damped wave equation
\begin{equation}\label{1}
\left\{\begin{array}{ll} \,\, \displaystyle {u_{tt}-\Delta
u -\frac{b_0}{(1+t)^{\beta}}\Delta u_t=|u|^p} &\displaystyle {x\in {\mathbb{R}^n},\,t>0,}\\
{}\\
\displaystyle{u(0,x)= u_0(x),\;\;u_t(0,x)=
u_1(x),\qquad\qquad}&\displaystyle{x\in {\mathbb{R}^n},}
\end{array}
\right.
\end{equation}
where $n\geq 1$, $p>1$, $b_0$  is a positive constant, and $\beta\in\mathbb{R}$. Without loss of generality, we assume that $b_0=1$.\\
 Throughout this paper, we assume that 
\begin{equation}\label{Existence}
\left\{\begin{array}{ll}
p\in(1,\infty)&\hbox{for}\,\, n=1,2,\\
{}\\
 p\in(1, \frac{n}{n-2}]&\hbox{for}\,\, n\geq 3,\\
\end{array}
\right.
\end{equation}
and the initial data are in the energy space
\begin{equation}\label{2}
(u_0,u_1)\in H^1(\mathbb{R}^n)\times L^2(\mathbb{R}^n).
\end{equation}
Here after, $\|\cdotp\|_q$ and $\|\cdotp\|_{H^1}$ $(1\leq q\leq
\infty)$ stand for the usual $L^q(\mathbb{R}^n)$-norm and
$H^1(\mathbb{R}^n)$-norm, respectively.

\medskip

In this paper, we study the blow-up result of  solution of \eqref{1}. Before
going on, it is necessary to mention that the
 case  $b_0=0$  in \eqref{1} is   the classical semilinear wave equation for which we have the Strauss conjecture. More precisely, this case  is characterized by a critical power, denoted by $p_S$,  which is the positive  solution of the following quadratic equation
$(n-1)p^2-(n+1)p-2=0,$
and is given  by
\begin{equation}
p_S=p_S(n):=\frac{n+1+\sqrt{n^2+10n-7}}{2(n-1)}.
\end{equation}
More precisely,  if $p \le p_S$ then there is no global solution for  \eqref{1} under suitable sign assumptions for the initial data, 
and for $p > p_S$ a global solution exists for small initial data; see e.g. \cite{John2,Strauss,YZ06,Zhou} among many other references. \\
A slightly less sharp blow-up result under much weaker assumptions was obtained by Kato \cite{Kato} with a much easier proof. In particular, Kato pointed out the role of the exponent ${(n+1)}/{(n-1)}<p_S(n),$ for $n\geq2,$ in order to have more general initial data, but still with compact support.\\
 We take the opportunity to mention here that the test function method, introduced by \cite{Zhang} and used by \cite{FinoKarch, Finokirane, PM}, plays a similar role as of Kato's method to prove blow-up results. In fact, the test function is effective in the case of parabolic equations which means that it provide us exactly the critical exponent $p_c$, while in the case of hyperbolic equations (cf. \cite{PM}) we get the so-called Kato's exponent $p^*$ i.e. we obtain a blow-up result for $p\leq p^*< p_c$. This is one of the weakness of the test function method but in general it can be applied to a more general equation and system.\\

When $\beta=0$, and $b_0=1$, problem \eqref{1} is reduced to
\begin{equation}\label{strongdamping}
 \displaystyle {u_{tt}-\Delta
u -\Delta u_t=|u|^p}, \ \ \ \ \displaystyle {x\in {\mathbb{R}^n},t>0,}
\end{equation}
which is called the viscoelastic damping case. D'Ambrosio and Lucente \cite[Theorem~4.2]{Ambrosio} proved that the solution of \eqref{strongdamping} blows-up in finite time when $1<p\leq (n+1)/(n-1)_+$, where $(\cdotp)_+:= \max\{0,\cdotp\}$, by applying the test function method. Similar result has been obtained recently  by Fino \cite{Fino} in the case of an exterior domain. On the other hand, D'Abbicco-Reissig \cite{Dabbicco} proved that there exists a global solution for \eqref{strongdamping} when $p>1+\frac{3}{n-1}$ $(n\geq 2)$ for sufficiently small initial data. Therefore, the exact value of the critical exponent is still an open question.
\\

When $\beta\neq0$, and $b_0=1$, we give an intuitive observation for understanding the influence of the damping  term ($-\frac{1}{{(1+t)}^{\beta}}\Delta u_t$) by scaling argument. Let $u(t,x)$ be a solution of the linear strong damped wave equation
\begin{equation}\label{linear}
u_{tt}(t,x)-\Delta
u(t,x) -\frac{1}{(1+t)^{\beta}}\Delta u_t(t,x)=0.
\end{equation}
When $\beta\geq-1$, we put
\begin{equation}\label{transf}
v(t,x)=u(\lambda(1+ t),\lambda x), \qquad \lambda (t+1)=s,\,\,\lambda x=y,
\end{equation}
%
with a parameter $\lambda>0$, we have
$$v_{ss}(s,y)-\Delta
v(s,y)-\frac1{\lambda^{\beta+1}s^{\beta}}\Delta v_s(s,y)=0.$$
Thus, when $\beta=-1$ we notice that the equation \eqref{linear} is invariant, while when $\beta>-1$, letting $\lambda\rightarrow \infty$, we obtain the wave equation without damping
$$v_{ss}(s,y)-\Delta v(s,y)=0.$$
We note that $\lambda\rightarrow \infty$ is corresponding to $t\rightarrow+\infty$.\\
On the other hand, when $\beta<-1$, we put
$$v(t,x)=u(\lambda^{\frac{2}{1-\beta}}(1+ t),\lambda x), \qquad \lambda ( t+1)=s,\,\,\lambda x=y,$$
%
with a parameter $\lambda>0$, we have
$$v_{ss}(s,y)-\frac{1}{s^{\beta}}\Delta v_s(s,y)-\lambda^{\frac{2(\beta+1)}{1-\beta}}\Delta
v(s,y)=0.$$
In this case, letting $\lambda\rightarrow \infty$, we obtain the hyperbolic equation
$$v_{ss}(s,y)-\frac{1}{s^{\beta}}\Delta v_s(s,y)=0.$$

In this paper, our goal is to generalize Kato's exponent and give sufficient conditions for finite time blow-up of a new type of class of equations \eqref{1} for $b_0> 0$, $\beta\in\mathbb{R}$. Let us mention that our blow-up results and initial conditions are similar to that of Kato.\\

This paper is organized as follows. We start  in Section \ref{sec-main} by introducing the mild solution of \eqref{1}.
Then, we state the main theorem of our work.  In Section \ref{loc},
we study the local existence of the solutions of equation \eqref{1} in the energy space $H^1(\mathbb{R}^n)\times L^2(\mathbb{R}^n).$ Finally, in Section \ref{bmain}, we prove the blow-up theorem (Theorem \ref{Blow-up}).

\section{Main results}\label{sec-main}
\par
This section is aimed to state our main results. For that purpose, we first start by giving the definition of the mild solution of \eqref{1}.
\begin{definition}(Mild solution)\\
Let $(u_0,u_1)\in H^1(\mathbb{R}^n)\times L^2(\mathbb{R}^n)$. We say that a function 
$$u\in C([0,T];H^1(\mathbb{R}^n))\cap C^1([0,T];L^2(\mathbb{R}^n))$$
is a mild solution of (\ref{1})  if and $u$ has the initial data $u(0)=u_0$, $u_t(0)=u_1$ and satisfies
the integral equation
 \begin{equation}\label{mild}
u(t,x)=R(t,0)(u_0,u_1)+\int_0^tS(t,s)|u(s)|^{p}\,ds
\end{equation}
in the sense of $H^1(\mathbb{R}^n)$, where the operators $R$ and $S$ are defined below. Moreover, if $T>0$ can be arbitrary chosen, then $u$ is called a global mild solution (\ref{1}).
\end{definition}

Here is the statement of our main theorem in this paper.

\begin{theorem}[\textbf{Blow-up}] \label{Blow-up}
 We assume that
$$(u_0,u_1)\in \big(L^1(\mathbb{R}^n)\cap H^1(\mathbb{R}^n)\big) \times \big(L^1(\mathbb{R}^n)\cap L^2(\mathbb{R}^n)\big) $$
satisfying the following condition:
\begin{equation} \label{initialdata}
\int_{\mathbb{R}^n} u_1(x)dx > 0.
\end{equation}
If 
\begin{equation}\label{30oc}
\left\{\begin{array}{ll}
 1<p\leq \frac{n+1}{(n-1)_+}&\text{if}\,~ \beta\geq-1,\\
 &\\
1<p\leq \frac{n(1-\beta)+2}{(n(1-\beta)-2)_+}&\text{if}\,~ \beta\leq-1,\\
\end{array}
\right.
\end{equation}
where $(\cdotp)_+:= \max\{0,\cdotp\}$, then the mild solution of \eqref{1} blows-up in finite time.
\end{theorem}

\begin{rmk}\label{rk1bis}
We   stress that 
  the exponent  $\frac{n+1}{(n-1)_+}$ appearing  in \eqref{30oc} was  introduced  first in \cite{Kato} 
  to prove the nonexistence of global
solutions to the semilinear wave equation with the nonlinearity $|u|^p,$ for small initial data with compact support. 
\end{rmk}

\begin{rmk}\label{rk1bis1}
Theorem 1 asserts that if $\beta \ge  -1 $, then  the critical exponent for $p$, is greater than or equal to $\frac{n+1}{(n-1)_+}$.
 Therefore, the   
blow-up region obtained in the present work in the case 
$\beta\ge -1$ constitute somehow an extension to the results related to the blow-up region of the solution of the equation \eqref{strongdamping}   obtained in \cite{Ambrosio}.  
\end{rmk}

\begin{rmk}\label{rk1}
 It is interesting to recall that thanks to the transformation  \eqref{transf} 
 the asymptotic behavior of the solution to \eqref{linear}  in the case $\beta \ge  -1 $, is given by the free wave equation.  
 Unfortunatly,  the  time-dependent damping  term ($-\frac{1}{{(1+t)}^{\beta}}\Delta u_t$) 
 makes the problem parabolic and loses the property that the speed of propagation is finite. 
 For this reason,  we obtain     the blow-up region    given by  \eqref{30oc} which is included  in $(1, p_S(n)]$.
 We think  it is likely possible to
  extend the blow-up region  to a larger interval $(1,p_c]$ where  
  $\frac{n+1}{(n-1)_+}<p_c\le p_S(n)$  at least for $\beta$ large enough.  
\end{rmk}

\section{Local existence}\label{loc}
To prove  that the Cauchy problem for  (\ref{1}) is locally well-posed in the space $H^1(\mathbb{R}^n)\times L^2(\mathbb{R}^n),$   it  is a natural to start by studying the linear homogeneous case
\subsection{Linear homogeneous case}\label{LHC}
We consider the linear homogeneous equation
\begin{equation}\label{eq2.1}
\left\{\begin{array}{ll}
\,\, \displaystyle {u_{tt}-\Delta
u -
\frac{1}{(1+t)^{\beta}}\Delta u_t =0,} &\displaystyle {t>0,x\in {\mathbb{R}^n},}\\
{}\\
\displaystyle{u(0,x)=  u_0(x),\;\;u_t(0,x)=  u_1(x),\qquad\qquad}&\displaystyle{x\in {\mathbb{R}^n},}\\
\end{array}
\right.
\end{equation}
\begin{definition}(Strong solution)\\
Let $(u_0,u_1)\in (H^2(\mathbb{R}^n))^2$. A function $u$ is said to be a strong solution of $(\ref{eq2.1})$ if
$$u\in C^{1}([0,\infty);H^2(\mathbb{R}^n))\cap C^{2}([0,\infty);L^2(\mathbb{R}^n)),$$
and $u$ has the initial data $u(0)=u_0$, $u_t(0)=u_1$ and satisfies
the equation (\ref{eq2.1}) in the sense of $L^2(\mathbb{R}^n)$.
\end{definition}
\begin{proposition}\label{prop2.1}
For each $(u_0,u_1)\in (H^2(\mathbb{R}^n))^2$, there exists a unique strong solution $u$ of problem $(\ref{eq2.1})$ that satisfies the following energy estimates
\begin{equation}\label{eq2.2}
\int_{\mathbb{R}^n}(u_t^2(t,x)+|\nabla u(t,x)|^2)\,dx\leq \int_{\mathbb{R}^n}(u_1^2(x)+|\nabla u_0(x)|^2)\,dx,
\end{equation}
\begin{equation}\label{eq2.3}
\|u(t)\|_{L^2}\leq \|u_0\|_{L^2}+T\|(u_1,\nabla
u_0)\|_{L^2\times
L^2},
\end{equation}
for any $T>0$, and all $0\leq t\leq T$.
\end{proposition}
\proof  $\,$ The existence of the strong solution can be done easily by the semigroup theory (cf. \cite{CH}).  
We now  focus  to the proof of estimates \eqref{eq2.2} and \eqref{eq2.3}. 
By multiplying   (\ref{eq2.1}) by $u_t$, integrating over $\mathbb{R}^n$ 
and   performing some  integration by parts in space, 
we obtain
\begin{equation}\label{A1}
\frac{1}{2}\frac{d}{dt}\int_{\mathbb{R}^n} (u^2_t+|\nabla u|^2)\,dx+
\frac{1}{(1+t)^{\beta}}\int_{\mathbb{R}^n} |\nabla u_t|^2\,dx=0,
\qquad   \forall
t\ge 0.
\end{equation}
By integrating in time between $0$ and $t$  the equality \eqref{A1}, we get
\begin{equation*}
\int_{\mathbb{R}^n} (u^2_t+|\nabla u|^2)\,dx+2\int_0^t
\frac{1}{(1+s)^{\beta}}\int_{\mathbb{R}^n} |\nabla u_s|^2\,dx\,ds
= \int_{\mathbb{R}^n} (u^2_1+|\nabla u_0|^2)\,dx, \qquad   
 \forall
t\ge 0,
\end{equation*}
which implies the estimate  \eqref{eq2.2}. Next, we prove (\ref{eq2.3}). 
Thanks to the basic identity $u(t)=u_0+\int_0^tu_s(s)\,ds,$
we conclude
\begin{equation}\label{A2}
\|u(t)\|_{L^2}\leq \|u_0\|_{L^2}+\int_0^t\|u_s(s)\|_{L^2}\,ds \qquad   
 \forall
t\ge 0.
\end{equation}
By combining   \eqref{eq2.2} and   \eqref{A2}, we infer
\begin{eqnarray*}
\|u(t)\|_{L^2}\leq  \|u_0\|_{L^2}+\int_0^t\|(u_1,\nabla
u_0)\|_{L^2\times
L^2}\,ds
\leq \|u_0\|_{L^2}+T\|(u_1,\nabla
u_0)\|_{L^2\times L^2} \qquad   
 \forall
t\ge 0.
\end{eqnarray*}
This follows \eqref{eq2.3} and we complete the proof of
Proposition \ref{prop2.1}.  $\hfill\blacksquare$\\

Let us denote by $R(t)$ the operator which maps the initial data $(u_0,u_1)\in
(H^2(\mathbb{R}^n))^2$ to the strong solution $u(t)\in
H^2(\mathbb{R}^n)$ at the time $t \geq 0$, i.e. the solution $u$ of (\ref{eq2.1})
is defined by $u(t)=R(t)(u_0,u_1)$. 
\begin{rmk}\label{rmk2.1}
From Proposition \ref{prop2.1}, the operator $R(t)$ can be
extended uniquely such that $R(t):H^1(\mathbb{R}^n)\times
L^2(\mathbb{R}^n)\longrightarrow C([0,\infty),H^1(\mathbb{R}^n))\cap
C^1([0,\infty),L^2(\mathbb{R}^n))$. Indeed, for any fixed $T>0$, due to the energy estimates (\ref{eq2.2})-(\ref{eq2.3}), the following estimation
$$
\|R(t)(u_0,u_1)\|_{H^1}+\|\partial_t(R(t)(u_0,u_1))\|_{L^2}\leq
C(1+T)\|(u_0,u_1)\|_{H^1\times L^2},
$$
holds for all $0\leq t\leq T$. It follows that the operator $R(t)$ can be extended uniquely to an operator such that $R(t):H^1(\mathbb{R}^n)\times
L^2(\mathbb{R}^n)\longrightarrow C([0,T],H^1(\mathbb{R}^n))\cap
C^1([0,T],L^2(\mathbb{R}^n))$. Since $T$ is arbitrary, we conclude the desired extension.
\end{rmk}

\subsection{Linear inhomogeneous case}
Let us now consider the linear inhomogeneous equation
\begin{equation}\label{eq3.3}
\left\{\begin{array}{ll} \,\, \displaystyle {u_{tt}-\Delta
u -
\frac{1}{(1+t)^{\beta}}\Delta u_t=F(t,x),} &\displaystyle {t>0,x\in {\mathbb{R}^n},}\\
{}\\
\displaystyle{u(0,x)= u_0(x),\;\;u_t(0,x)=
u_1(x),\qquad\qquad}&\displaystyle{x\in {\mathbb{R}^n}.}
\end{array}
\right.
\end{equation}

\begin{definition}
Let $(u_0,u_1)\in H^1\times L^2$ and $F\in C([0,\infty);L^2)$. We
say that a function $u$ is a mild solution of (\ref{eq3.3}) if $u\in
C([0,\infty);H^1)\cap C^1([0,\infty);L^2)$ and $u$ has the initial
data $u(0)=u_0$, $u_t(0)=u_1$ and satisfies the integral equation
 \begin{equation}\label{eq3.6}
u(t,x)=R(t,0)(u_0,u_1)+\int_0^tS(t,s)F(s,x)\,ds
\end{equation}
in the sense of $H^1(\mathbb{R}^n)$, where $S(t,s)g:=R(t,s)(0,g)$
for a function $g\in H^1(\mathbb{R}^n)$.
\end{definition}
By a classical result as in \cite{CH, Pazy} or similarly as in
 \cite[Proposition 9.15]{Yuta}, we have the following
\begin{proposition}\label{eq3.5}${}$\\
Let $(u_0,u_1)\in H^1(\mathbb{R}^n)\times L^2(\mathbb{R}^n)$, $F\in C([0,\infty);L^2(\mathbb{R}^n))$. Then
there exists a unique mild solution $u$ of (\ref{eq3.3}). Moreover,
the mild solution $u$ satisfies the following energy estimates 
\begin{equation}\label{}
\|(u_t,\nabla u)(t)\|_{L^2\times L^2}\leq C\|(u_1,\nabla
u_0)\|_{L^2\times L^2}+C\int_0^t\|F(s,\cdot)\|_{L^2}\,ds,
\end{equation}
\begin{equation}\label{}
\|u(t)\|_{L^2}\leq C\|u_0\|_{L^2}+C\,t\,\|(u_1,\nabla
u_0)\|_{L^2\times
L^2}+\int_0^t\int_0^s\|F(\tau,\cdot)\|_{L^2}\,d\tau\,ds.
\end{equation}

\end{proposition}

\subsection{Semilinear case}

Using Gagliardo-Nirenberg's inequality, Proposition \ref{eq3.5} and the Banach fixed point theorem we get the following local existence theorem.
\begin{proposition}\label{prop1}
Let $\beta\in\mathbb{R}$. Under the assumptions
\eqref{Existence}-\eqref{2}, the problem $(\ref{1})$
admits a unique maximal mild solution $u$, i.e. satisfies the
integral equation (\ref{mild}) such that
$$u\in C([0,T_{\max}),H^1(\mathbb{R}^n))\cap C^1([0,T_{\max}), L^2(\mathbb{R}^n)),$$
where $0< T_{\max}\leq\infty$. Moreover, if $T_{\max}<\infty$, then
it follows that
$$\|u(t)\|_{H^1}+\|u_t(t)\|_2\rightarrow\infty\qquad\mbox{as}\;\;
t\rightarrow T_{\max}.$$
\end{proposition}

\section{Proof of Theorem \ref{Blow-up}}\label{bmain}
In order to prove Theorem \ref{Blow-up}, we are going to use the test function method which is rely on the weak solution of (\ref{1}). More precisely, the weak formulation associated with  (\ref{1}) reads  as follows:

\begin{definition}(Weak solution)\\
Let $T>0$, and $u_0,u_1\in L_{loc}^1(\mathbb{R}^n)$. A function $u$ is said to be a weak solution of $(\ref{1})$ if
$$u\in L^p((0,T);L_{loc}^p(\mathbb{R}^n)),$$
 and $u$ satisfies the weak formulation
 \begin{eqnarray*}
&{}&\int_0^T\int_{\mathbb{R}^n}|u|^p\psi\,dx\,dt+\int_{\mathbb{R}^n}u_1(x)\psi(0,x)\,dx-\int_{\mathbb{R}^n}u_0(x)\Delta\psi(0,x)\,dx-\int_{\mathbb{R}^n}u_0(x)\psi_t(0,x)\,dx\\
&{}&=\int_0^T\int_{\mathbb{R}^n}u\psi_{tt}\,dx\,dt+\int_0^T\frac{1}{(1+t)^{\beta}}\int_{\mathbb{R}^n}u\,\Delta\psi_t\,dx\,dt-\int_0^T\int_{\mathbb{R}^n}u \Delta \psi\,dx\,dt-\int_0^T\frac{\beta}{(1+t)^{\beta+1}}\int_{\mathbb{R}^n}u\,\Delta \psi\,dx\,dt,
\end{eqnarray*}
for all compactly supported function $\psi\in
C^2([0,T]\times\mathbb{R}^n)$ such that $\psi(\cdotp,T)=0$ and
$\psi_t(\cdotp,T)=0$. We denote the lifespan for the weak solution by
$$T_w(u_0):=\sup\{T\in(0,\infty];\,\,\hbox{there exists a unique weak solution u to (\ref{1})}\}.$$
Moreover, if $T>0$ can be arbitrary chosen, i.e. $T_w(u_0)=\infty$, then $u$ is called a global weak solution of (\ref{1}).
\end{definition}
We also need the following 
\begin{rmk}\label{mildweak}$(\mbox{Mild $\rightarrow$ Weak})\;$\\
 Let $(u_0,u_1)\in H^1(\mathbb{R}^n)\times L^2(\mathbb{R}^n)$. Under the assumption $(\ref{Existence})$, if $u$ is a global  mild solution of $(\ref{1})$, then $u$ is a  global  weak solution of \eqref{1}.
\end{rmk}

\noindent{\bf Proof of Theorem \ref{Blow-up}.}
Let $u$ a global  mild solution of \eqref{1}.Thanks to  Remark \ref{mildweak},  we conclude that $u$ is a  global weak solution of \eqref{1}. i.e.
  \begin{eqnarray}\label{weak1}
&{}&\int_0^T\int_{\mathbb{R}^n}|u|^p\psi\,dx\,dt+\int_{\mathbb{R}^n}u_1(x)\psi(0,x)\,dx-\int_{\mathbb{R}^n}u_0(x)\Delta\psi(0,x)\,dx-\int_{\mathbb{R}^n}u_0(x)\psi_t(0,x)\,dx\\
&{}&=\int_0^T\int_{\mathbb{R}^n}u\psi_{tt}\,dx\,dt+\int_0^T
\frac{1}{(1+t)^{\beta}}\int_{\mathbb{R}^n}u\,\Delta\psi_t\,dx\,dt-\int_0^T\int_{\mathbb{R}^n}u \Delta \psi\,dx\,dt-\int_0^T
\frac{\beta}{(1+t)^{\beta+1}}\int_{\mathbb{R}^n}u\,\Delta \psi\,dx\,dt,\nonumber
\end{eqnarray}
for all $T>0$, and all compactly supported function $\psi\in
C^2([0,T]\times\mathbb{R}^n)$ such that $\psi(\cdotp,T)=0$ and
$\psi_t(\cdotp,T)=0$.

Let $T>0$. Now, we introduce the following test function:
\begin{equation}\label{psi}
\psi(x,t)=\psi^\ell_1(x)\psi^\eta_2(t)
\end{equation}
where
\begin{equation}\label{psi1}
\psi_1(x):=\Phi\left(\frac{|x|}{T^d}\right),\qquad \psi_2(t):=\Phi\left(\frac{t}{T}\right),
\end{equation}
where  $\ell,\eta$ are sufficiently large constants that will be determined later
 and $\Phi\in C^\infty(\mathbb{R}_+)$ be a cut-off non-increasing
function such that
$$\Phi(r)=\left\{\begin {array}{ll}\displaystyle{1}&\displaystyle{\quad\mbox{if }\,0\leq r\leq 1/2,}\\\\
\displaystyle{\searrow}&\displaystyle{\quad\mbox{if }\,1/2\leq r\leq 1,}\\\\
\displaystyle{0}&\displaystyle{\quad\mbox {if }\,r\geq 1.}
\end {array}\right.$$\\
The
constant $d>0$ in the definition of $\psi_1$ is fixed and will be
chosen later. In the following, we denote by $\Omega(T)$ the support
of $\psi_1$ and by $\Delta(T)$ the set containing the support of
$\Delta\psi_1$ which are defined as follows:
$$\Omega(T)=\{x\in\mathbb{R}^n:\;|x|\leq 2T^d\},\quad \Delta(T)=\{x\in\mathbb{R}^n:\;T^d/2\leq|x|\leq T^d\}.$$
By \eqref{weak1}, we  get that
\begin{eqnarray}\label{3}
&{}&\int_0^T\int_{\Omega(T)}|u|^p\psi\,dx\,dt+\int_{\Omega(T)}u_1(x)\psi(0,x)\,dx\nonumber\\
&{}&\leq \int_{\frac{T}{2}}^T\int_{\Omega(T)}|u|\,|\psi_{tt}|\,dx\,dt+\int_0^T\int_{\Delta(T)}|u|\,|\Delta \psi|\,dx\,dt+\int_{\frac{T}{2}}^T
\frac{1}{(1+t)^{\beta}}\int_{\Delta(T)}|u||\Delta\psi_t|\,dx\,dt\nonumber\\
&{}&\,\,\,\,\,\,+\int_0^T
\frac{\beta}{(1+t)^{\beta+1}}\int_{\Delta(T)}|u|\,\,|\Delta \psi|\,dx\,dt+\int_{\Omega(T)}|u_0|\left(|\Delta\psi(0,x)|+|\psi_t(0,x)|\right)\,dx\nonumber\\
&{}&=:I_1+I_2+I_3+I_4+I_5.
\end{eqnarray}
Let $\varepsilon>0$. By applying $\varepsilon$-Young's inequality 
$$
AB\leq\varepsilon A^p+C(\varepsilon,p)B^{p^\prime},\quad A\geq0,\;B\geq0,\;p+p^\prime=pp^\prime,\quad C(\varepsilon,p)=\varepsilon^{-1/(p-1)} p^{-p/(p-1)} (p-1),
$$
we obtain
\begin{eqnarray}\label{4}
I_1&\leq&\int_0^T\int_{\Omega(T)}|u|\psi^{1/p}\psi^{-1/p}\psi_1^\ell
|(\psi_2^\eta)_{tt}|\,dx\,dt\nonumber\\
&\leq&\varepsilon\int_0^T\int_{\Omega(T)}|u|^{p}\psi\,dx\,dt+C\int_0^T\int_{\Omega(T)}\psi_1^\ell\psi_2^{-\eta/(p-1)}
|(\psi_2^\eta)_{tt}|^{p'}\,dx\,dt.
\end{eqnarray} 
Using 
$(\psi_2^\eta)_{tt}=\eta\psi_2^{\eta-1}(\psi_2)_{tt}+\eta(\eta-1)\psi_2^{\eta-2}|(\psi_2)_t|^2,$  the inequality \eqref{4} becomes
\begin{equation}\label{4bis}
I_1
\leq\varepsilon\int_0^T\int_{\Omega(T)}|u|^{p}\psi\,dx\,dt+C\int_0^T\int_{\Omega(T)}\psi_1^\ell\psi_2^{\eta-p'}|(\psi_2)_{tt}|^{p'}\,dx\,dt+C\int_0^T\int_{\Omega(T)}\psi_1^\ell\psi_2^{\eta-2p'}|(\psi_2)_{t}|^{2p'}\,dx\,dt.
\end{equation} 
 Using and  proceeding similarly as for \eqref{4bis} and using the identity   $\Delta(\psi_1^\ell)=\ell\psi_1^{\ell-1}\Delta\psi_1+\ell(\ell-1)\psi_1^{\ell-2}|\nabla\psi_1|^2,$  we easily deduce
\begin{equation}\label{5}
I_2
\leq\varepsilon\int_0^T\int_{\Omega(T)}|u|^{p}\psi\,dx\,dt+C\int_0^T\int_{\Omega(T)}\psi_2^\eta\psi_1^{\ell-p'}|\Delta\psi_1|^{p'}\,dx\,dt+C\int_0^T\int_{\Omega(T)}\psi_2^\eta\psi_1^{\ell-2p'}|\nabla\psi_1|^{2p'}\,dx\,dt.
\end{equation} 
In the same way, thanks to $\Delta\psi_t= \Delta (\psi_1^\ell)(\psi^{\eta}_2)_t$, we write
\begin{eqnarray}\label{6}
I_3&\leq&\int_{\frac{T}{2}}^T
\frac{1}{(1+t)^{\beta}}\int_{\Omega(T)}|u|\psi^{1/p}\psi^{-1/p}\,|\Delta(\psi_1^\ell)|\,|(\psi^\eta_2)_t|\,dx\,dt\nonumber\\
&\leq&\varepsilon\int_0^T\int_{\Omega(T)}|u|^{p}\psi\,dx\,dt+C\int_{\frac{T}{2}}^T
\frac{1}{(1+t)^{\beta p'}}\int_{\Omega(T)}
\psi_1^{-\ell/(p-1)}
|\Delta(\psi_1^\ell)|^{p'}\,\psi_2^{-\eta/(p-1)}\,|(\psi^\eta_2)_t|^{p'}\,dx\,dt.
\end{eqnarray} 
 Clearly,
$$
\frac{1}{(1+t)^{\beta p^\prime}}\leq C\,T^{-\beta p'},\qquad\forall\,t\in\left(\frac{T}{2},T\right).$$
 Therefore, 
\begin{eqnarray}\label{6bis}
I_3
&\leq&\varepsilon\int_0^T\int_{\Omega(T)}|u|^{p}\psi\,dx\,dt+C\,T^{-\beta p'}\,\int_0^T\int_{\Omega(T)}\psi_2^{\eta-p'}|(\psi_2)_t|^{p'}\psi_1^{\ell-p'}|\Delta\psi_1|^{p'}\,dx\,dt\nonumber\\
&{}&\,\,+\,C\,T^{-\beta p'}\,\int_0^T\int_{\Omega(T)}\psi_2^{\eta-p'}|(\psi_2)_t|^{p'}\psi_1^{\ell-2p'}|\nabla\psi_1|^{2p'}\,dx\,dt.
\end{eqnarray} 
In the same manner,
\begin{eqnarray}\label{7}
I_4&\leq&C\int_{0}^T\frac{1}{(1+t)^{\beta+1}}\int_{\Omega(T)}|u|\psi^{1/p}\psi^{-1/p}\,\,\psi^\eta_2\,|\Delta(\psi_1^\ell)|\,dx\,dt\nonumber\\
&\leq&\varepsilon\int_0^T\int_{\Omega(T)}|u|^{p}\psi\,dx\,dt+C\,\int_0^T\frac{1}{(1+t)^{(\beta+1)p'}}\int_{\Omega(T)}
\psi_2^\eta\psi_1^{\ell-p'}|\Delta \psi_1|^{p'}\,dx\,dt\nonumber\\
&{}&\,\,+\,C\,\int_0^T\frac{1}{(1+t)^{(\beta+1)p'}}\int_{\Omega(T)}\psi_2^{\eta}\psi_1^{\ell-2p'}|\nabla \psi_1|^{2p'}\,dx\,dt.
\end{eqnarray} 
Finally, it remains only to control the term $I_5$. Note 
by exploiting  the identities 
$$\Delta \psi (0,x)=\Delta(\psi_1^\ell)=\ell\psi_1^{\ell-1}\Delta\psi_1+\ell(\ell-1)\psi_1^{\ell-2}|\nabla\psi_1|^2\quad\hbox{and}\quad \psi_t (0,x)=\eta \psi_1^\ell(\psi_2)_t(0),$$
 we infer
\begin{equation}\label{8}
I_5 \leq C\,\int_{\Omega(T)}|u_0|\left(\psi_1^{\ell-1}|\Delta \psi_1|+\psi_1^{\ell-2}|
\nabla \psi_1|^2+|(\psi_2)_t(0)|\psi_1^\ell\right)\,dx.
\end{equation} 
Plugging \eqref{3} together  with \eqref{4bis}-\eqref{8} and  choosing $\varepsilon$ small enough, we deduce  that
\begin{eqnarray}\label{9}
&{}&\int_0^T\int_{\Omega(T)}|u|^p\psi\,dx\,dt+\int_{\Omega(T)}u_1(x)\psi(0,x)\,dx\nonumber\\
&{}&\leq C\int_0^T\int_{\Omega(T)}\psi_1^\ell\psi_2^{\eta-p'}|(\psi_2)_{tt}|^{p'}\,dx\,dt+C\int_0^T\int_{\Omega(T)}\psi_1^\ell\psi_2^{\eta-2p'}|(\psi_2)_{t}|^{2p'}\,dx\,dt\nonumber\\
&{}&\quad+\,C\int_0^T\int_{\Omega(T)}\psi_2^\eta\psi_1^{\ell-p'}|\Delta\psi_1|^{p'}\,dx\,dt+C\int_0^T\int_{\Omega(T)}\psi_2^\eta\psi_1^{\ell-2p'}|\nabla\psi_1|^{2p'}\,dx\,dt\nonumber\\
&{}&\quad+\,C\,T^{-\beta p'}\,\int_0^T\int_{\Omega(T)}\psi_2^{\eta-p'}|(\psi_2)_t|^{p'}\psi_1^{\ell-p'}|\Delta\psi_1|^{p'}\,dx\,dt\nonumber\\
&{}&\quad
+\,C\,T^{-\beta p'}\,\int_0^T\int_{\Omega(T)}\psi_2^{\eta-p'}|(\psi_2)_t|^{p'}\psi_1^{\ell-2p'}|\nabla\psi_1|^{2p'}\,dx\,dt\nonumber\\
&{}&\quad
+\,C\,\int_0^T\frac{1}{(1+t)^{\frac{(\beta+1)p}{p-1}}}\int_{\Omega(T)}
\psi_2^\eta\psi_1^{\ell-p'}|\Delta \psi_1|^{p'}\,dx\,dt
+C\,\int_0^T\frac{1}{(1+t)^{\frac{(\beta+1)p}{p-1}}}\int_{\Omega(T)}\psi_2^{\eta}\psi_1^{\ell-2p'}|\nabla \psi_1|^{2p'}\,dx\,dt\nonumber\\
&{}&\quad+\,C\,\int_{\Omega(T)}|u_0|\left(\psi_1^{\ell-1}|\Delta \psi_1|+\psi_1^{\ell-2}|\nabla \psi_1|^2+|(\psi_2)_t(0)|\psi_1^\ell\right)\,dx.
\end{eqnarray}
By neglecting the first term in the left-hand side and 
taking into account the expresion of $\psi$ given by \eqref{psi},
we obtain
\begin{eqnarray}\label{10}
\int_0^T\int_{\Omega(T)}|u|^p\psi\,dx\,dt+\int_{\Omega(T)}u_1(x)\psi^\ell_1(x)\,dx
&\leq& C\;T^{-2p'+1+nd}+C\;T^{-2dp'+1+nd}+C\;T^{-\beta p'-p'-2dp'+1+nd}\\
&& +\,C\;T^{-2dp'+nd}
\int_0^T(1+t)^{-\frac{(\beta+1)\,p}{p-1}}\,dt+\,C\left(T^{-2d}+T^{-1}\right)\,\int_{\mathbb{R}^n}|u_0(x)|\,dx.\nonumber
\end{eqnarray}

Now, we distinguish two cases:\\
 \noindent {\bf I. First case: $\beta\geq -1$.}\\
 In this case, we choose $d=1$. \\
  
\noindent \underline{Subcritical case $(p<\frac{n+1}{(n-1)_+})$.}\\
 Note that, as
\begin{equation}\label{11}
\int_0^T(1+t)^{-\frac{(\beta+1)\,p}{p-1}}\,dt\leq C\,\left\{\begin{array}{ll}
T^{1-\frac{(\beta+1)\,p}{p-1}}&\text{if}~\beta\,p<-1,\\
\ln(T)&\text{if}~\beta\,p=-1,\\
1&\text{if}~\beta\,p>-1,\\
\end{array}
\right.
\end{equation}
we have $\int_0^T(1+t)^{-\frac{(\beta+1)\,p}{p-1}}\,dt\leq C T$, for all $T>1$. Then 
  (\ref{10}) implies
\begin{equation}\label{10d1}
\int_0^T\int_{\Omega(T)}|u|^p\psi\,dx\,dt+\int_{\Omega(T)}u_1(x)\psi^\ell_1(x)\,dx
\leq C\;T^{-2p'+1+n}+C\;T^{-(\beta+1) p'-2p'+1+n}+
\,C\left(T^{-2}+T^{-1}\right)\,\int_{\mathbb{R}^n}|u_0(x)|\,dx,
\end{equation}
 for all $T>1$. We use the fact $\beta\ge -1$, to conclude that
\begin{equation}\label{10d11}
\int_0^T\int_{\Omega(T)}|u|^p\psi\,dx\,dt+\int_{\Omega(T)}u_1(x)\psi^\ell_1(x)\,dx
\leq C\;T^{-2p'+1+n}+
\,C\left(T^{-2}+T^{-1}\right)\,\int_{\mathbb{R}^n}|u_0(x)|\,dx, \quad \forall\,T>1.
\end{equation}

Note that,  we can easily see that
$-2p'+1+n<0,$
if
$p<\frac{n+1}{(n-1)_+}.$
Letting $T\rightarrow\infty$, and using the Lebesgue
dominated convergence theorem, we conclude that
$$\int_{\mathbb{R}^n}u_1(x)\,dx\leq 0.$$
This contradicts our assumption \eqref{initialdata}.\\

 \noindent \underline{Critical case $(p=\frac{n+1}{n-1})$, when $n\geq2$.}\\
 Let $n\geq2$, and 
 $$p=\frac{n+1}{n-1}.$$
 From the subcritical case \eqref{10d11}, we can see easily that we have
 \begin{equation}\label{regularity1}
 u\in L^p((0,\infty);L^p(\mathbb{R}^n)).
 \end{equation}
 On the other hand, by applying H\"{o}lder's inequality instead of Young's inequality, we get
$$\int_{\Omega(T)}u_1(x)\psi^\ell_1(x)\,dx \leq C\int_{\frac{T}{2}}^T\int_{\Omega(T^d)}|u|^{p}\psi\,dx\,dt+\,C\int_0^T\int_{\Delta(T^d)}|u|^{p}\psi\,dx\,dt+\,C\left(T^{-2}+T^{-1}\right)\,\int_{\mathbb{R}^n}|u_0(x)|\,dx.
$$
 Letting $T\longrightarrow\infty$ and taking into consideration \eqref{regularity1} we get a contradiction.\\

 \noindent {\bf II. Case of $\beta< -1$.}\\
In this case, we choose $d=\frac{1-\beta}{2}$. \\

\noindent Note that, as $p\beta<\beta<-1$, we have
\begin{equation}\label{bq2}
\int_0^T(1+t)^{-\frac{(\beta+1)\,p}{p-1}}\,dt\le C T^{1-{(\beta+1)\,p}}.
\end{equation}
Therfore, \eqref{10} becomes
\begin{eqnarray}\label{10b2}
\int_0^T\int_{\Omega(T)}|u|^p\psi\,dx\,dt+\int_{\Omega(T)}u_1(x)\psi^\ell_1(x)\,dx
&\leq& C\;T^{-2p'+1+nd}+C\;T^{-2dp'+1+nd}\nonumber\\
&{}&\, +\,C\;T^{-2dp'+nd+1-(\beta+1)\,p'}+\,C\left(T^{-2d}+T^{-1}\right)\,\int_{\mathbb{R}^n}|u_0(x)|\,dx.\nonumber
\end{eqnarray}

Now, we distinguish two subcases:\\

\noindent \underline{Subcritical case $p<\frac{n(1-\beta)+2}{n(1-\beta)-2}$}:\\

As $\beta<-1$, we have
$$-2dp'+1+nd<-(\beta+1) p'-2dp'+1+nd<0\quad\hbox{and}\quad -2p'+1+nd<0,,$$
where we have used
$$p<\frac{n(1-\beta)+2}{n(1-\beta)-2},$$
for all $n\geq1$. So, letting $T\longrightarrow\infty$, we get a contradiction.\\

\noindent \underline{Critical case $p=\frac{n(1-\beta)+2}{n(1-\beta)-2}$}:\\
 Let $n\geq1$, and 
 \begin{equation}\label{crit}
 p=\frac{n(1-\beta)+2}{n(1-\beta)-2}.
 \end{equation}
 From the subcritical case, we can see easily that we have
 \begin{equation}\label{regularity}
 u\in L^p((0,\infty);L^p(\mathbb{R}^n)).
 \end{equation}
 On the other hand, by applying H\"{o}lder's inequality instead of Young's inequality, we get
$$\int_{\Omega(T)}u_1(x)\psi(0,x)\,dx \leq C\int_{\frac{T}{2}}^T\int_{\Omega(T^d)}|u|^{p}\psi\,dx\,dt+\,C\int_0^T\int_{\Delta(T^d)}|u|^{p}\psi\,dx\,dt+\,C\left(T^{-2}+T^{-1}\right)\,\int_{\mathbb{R}^n}|u_0(x)|\,dx.
$$
 Letting $T\longrightarrow\infty$ and taking into consideration \eqref{regularity} we get a contradiction. $\hfill\blacksquare$\\

\section*{References}

\bibitem {CH} T. Cazenave, A. Haraux, \textit{Introduction aux probl\`emes
d'\'evolution semi-lin\'eaires}, Ellipses, Paris, (1990).

\bibitem {Dabbicco} {M. D'Abbicco, M. Reissig},  \textit{Semilinear structural damped waves}, Math. Meth. Appl. Sci. $\textbf{37}$ $(2014)$ $1570-1592$.

\bibitem {Ambrosio} {L. D'Ambrosio, S. Lucente},  \textit{Nonlinear Liouville theorems for Grushin and Tricomi operators}, Journal Differential Equations $\textbf{123}$ $(2003)$ 511-541.

\bibitem{Fino} {A. Z. Fino}, \textit{Finite time blow up for wave equations with strong damping in an exterior domain}, Mediterranean Journal of Mathematics $\textbf{17}$ $(2020)$, no. $6$, Paper No. $174$, $21$ pp.

\bibitem{FinoKarch} { A. Fino, G. Karch},  \textit{Decay of mass for nonlinear equation with fractional Laplacian},  J. Monatsh. Math. ${\textbf160}$ $(2010),$ $375-384.$
 
 \bibitem{Finokirane} {A. Z. Fino, M. Kirane},  \textit{Qualitative properties of solutions to a time-space fractional evolution equation}, J. Quarterly of Applied Mathematics $\textbf{70}$ $(2012)$, $133-157$.
 
\bibitem{John2}
{F. John}, {\it Blow-up of solutions of nonlinear wave equations in three space dimensions}, Manuscripta Math. {\bf 28} (1979), no. 1-3, 235--268.

\bibitem{Kato}{T. Kato}, \textit{Blow-up of solutions of some nonlinear hyperbolic equations}, Comp. Pure Appl. Math. ${\bf 33}$ $(1980),$ $501-505.$

 \bibitem{PM}  { E. Mitidieri, S. I. Pohozaev},   \textit{A priori estimates and
blow-up of solutions to nonlinear partial differential equations and
inequalities}, Proc. Steklov. Inst. Math. ${\textbf234}$ $(2001),$ $1-383.$

\bibitem{Strauss}
{W. A. Strauss}, {\it  Nonlinear scattering theory at low energy.} J. Functional Analysis, {\bf  41} (1981), no. 1, 110--133.

\bibitem{Yuta} Y. Wakasugi, \textit{On the diffusive structure for the damped wave
equation with variable coefficients}, Doctoral thesis, Osaka

\bibitem{YZ06}
{B. Yordanov and Q. S. Zhang}, {\it Finite time blow up for critical wave equations in high dimensions}, J. Funct. Anal., {\bf 231} (2006), 361--374.

\bibitem{Zhang} {Qi S. Zhang},  \textit{A blow up result for a nonlinear wave equation with damping:
the critical case}, C. R. Acad. Sci. Paris, Vol. $\textbf{333}$  $(2001),$ no. $2,$ $109-114.$

\bibitem{Zhou}
{Y. Zhou}, {\it  Blow up of solutions to semilinear wave equations with critical exponent in high dimensions}, Chin. Ann. Math. Ser. B {\bf 28} (2007), no. 2, 205--212.

\bibitem{Pazy}
{A.Pazy}, {\it Semi-groups of linear operators and applications to partial differential
equations.}  Appl. Math. Sci 44. Springer. New-York (1983).

\end{document}